\documentclass[12pt]{article}
\usepackage[utf8]{inputenc}

\usepackage{amssymb,xcolor,xspace}
\usepackage{amsmath,graphicx,comment}
\usepackage{tikz}
\usepackage{hyperref}

\usepackage{geometry} % to change the page dimensions
\usepackage{fancyhdr} % This should be set AFTER setting up the page geometry
\usepackage{sectsty}
\allsectionsfont{\sffamily\mdseries\upshape} % (See the fntguide.pdf for font help)
\usepackage[nottoc,notlof,notlot]{tocbibind} % Put the bibliography in the ToC
\usepackage[titles]{tocloft} % Alter the style of the Table of Contents

\def\cqfd{\skip10=\parfillskip\parfillskip=0pt
\enspace\hfill\symbolecqfd\par\parfillskip=\skip10\par\medskip}
\def\symbolecqfd{\rlap{$\sqcap$}$\sqcup$}

\newtheorem{theorem}{Theorem}[section]
\newtheorem{proposition}[theorem]{Proposition}

\newtheorem{corollary}[theorem]{Corollary}

\newtheorem{pro-fact}[theorem]{Fact}

\newtheorem{pro-example}[theorem]{Example}
\newenvironment{example}{\begin{pro-example}\rm}{\cqfd\end{pro-example}}

\newtheorem{pro-remark}[theorem]{Remark}
\newenvironment{remark}{\begin{pro-remark}\rm}{\cqfd\end{pro-remark}}

\newenvironment{preuve}{\rm \trivlist \item[\hskip \labelsep{\bf
Proof.}]}{\cqfd\endtrivlist}

\def\cqfd{\skip10=\parfillskip\parfillskip=0pt
\enspace\hfill\symbolecqfd\par\parfillskip=\skip10\par\medskip}
\def\symbolecqfd{\rlap{$\sqcap$}$\sqcup$}
\def\proof{\begin{preuve}}
\def\eop{\end{preuve}}

\let\phi\varphi
\def\llangle{\langle\!\langle}
\def\rrangle{\rangle\!\rangle}

\def\lcp{\mathop{\textsf{lcp}}}

\def\inv{^{-1}}
\let\epsilon\varepsilon
\def \calA {\mathcal{A}}
\def\red {\mathcal{R}}
\def\cred {\mathcal{CR}}
\def \O {\mathcal{O}}
\def \calP {\mathcal{P}}

\def\N{\mathbb{N}}
\def \P {\mathbb{P}}
\def\Nil{\textsf{N}}
\def\PSL{\textsf{PSL}}
\def \Q {\mathbb{Q}}
\def \R {\mathbb{R}}
\def \T {\mathbb{T}}
\def\Z{\mathbb{Z}}

\def\size{\textsf{size}}

\def\Aut{\textsf{Aut}}
\def\cogrowth{\textsf{cogrowth}}

\title{Random presentations and random subgroups: a survey\thanks{%
Part of this work was done while the third author was Ada Peluso Visiting Professor, Mathematics Department, Hunter College, CUNY. He also received support from the ANR (project \textsc{DeLTA}, ANR-16-CE40-0007).}}

\author{
    Fr\'ed\'erique Bassino, \small{\url{bassino@lipn.univ-paris13.fr}}\\
    \small{Univ. Paris 13, Sorbonne Paris Cité, LIPN, CNRS UMR7030, F-93430 Villetaneuse, France}%     
    \and
    Cyril Nicaud, \small{\url{nicaud@univ-mlv.fr}}\\
    \small{Université Paris-Est, LIGM (UMR 8049), CNRS, ENPC, ESIEE Paris, UPEM,}\\
    \small{F-77454, Marne-la-Vallée, France}%   
    \and
    Pascal Weil, \small{\url{pascal.weil@labri.fr}}\\
    \small{Univ. Bordeaux, LaBRI, CNRS UMR 5800, F-33400 Talence, France}\thanks{%
    LaBRI, Univ. Bordeaux, 351 cours de la Lib\'eration, 33400 Talence, France.}
    }

\begin{document}

\date{}
\maketitle
\tableofcontents

%%%%%%%%%%%%%%%%
\section{Introduction}

In infinite group theory, it is a classical and natural question to ask what most groups look like, what a random group looks like. The question can and must be made more precise: it is actually a question about random finitely presented groups, and in most of the literature, in fact a question about random finite group presentations on a fixed set of generators. The specific questions may be whether a random finite group presentation satisfies a small cancellation property, whether the group it presents is hyperbolic, residually finite, etc.

Early on, Gromov gave an answer to this question: almost all groups are hyperbolic (see \cite{1993:Gromov}, and \cite{1992:Olshanskii,1995:Champetier,2004:Ollivier} for precise statements and complete proofs).

When a group $G$ is fixed (\textit{e.g.}, the free group $F(A)$ over a given finite set $A$ of generators, a hyperbolic group, a braid group, the modular group), one may also ask what a random finitely generated subgroup looks like: is it free? is it malnormal? does it have finite index? in the case where $G = F(A)$, is the subgroup Whitehead minimal?

These questions have been abundantly studied in the literature. This paper is a partial survey and as such, it contains no new results, but it offers a synthetic view of a part of this very active field of research. We refer the reader to the survey by Ollivier \cite{2005:Ollivier} for more details on some of the topics discussed here, and to the survey by Dixon \cite{2002:Dixon} for a discussion of probabilistic methods in finite group theory.

A specific aspect of the present survey is that we discuss both random presentations and random subgroups, unlike Ollivier \cite{2005:Ollivier}. 

Random presentations were considered first in the literature, and we will start with them as well (Section~\ref{sec: presentations}). We then proceed to a discussion of results on random subgroups (Section~\ref{sec: random subgroups}). Finally, Section~\ref{sec: non-uniform distributions} discusses recent results on non-uniform distributions.

This is an updated and extended version of a talk given at the conference GAGTA-8 in Newcastle, NSW, in July 2014.

%%%%%%%%%%%%%%%%%%%
\subsection{Discrete representations}

The very notion of randomness relies on a notion of probability, and in many cases, on a notion of counting discrete representations of finitely presented groups, or finitely generated subgroups, of a certain size: how many subgroups of $F(A)$ are there, with a tuple of $f(n)$ generators of length at most $n$ for a given function $f$? how many whose Stallings graph (see Section~\ref{sec: stallings}) has at most $n$ vertices? how many isomorphism classes of 1-relator groups are there, whose relator has length at most $n$? etc\dots

So we must first discuss the discrete representations we will use to describe subgroups and presentations.

Let $A$ be a finite non-empty set and let $F(A)$ be the \emph{free group on $A$}. The symmetrized alphabet $\tilde A$ is the union of $A$ and a copy of $A$, $\bar A = \{\bar a \mid a\in A\}$, disjoint from $A$. We denote by $\tilde A^*$ the set of all words on the alphabet $\tilde A$. The operation $x\mapsto \bar x$ is extended to $\tilde A^*$ by letting $\bar{\bar a} = a$ and $\overline{ua} = \bar a\bar u$ for all $a\in \tilde A$ and $u\in \tilde A^*$. Recall that a word is \emph{reduced} if it does not contain a factor of the form $a\bar a$ ($a\in \tilde A$). The (free group) \emph{reduction} of a word $u \in \tilde A^*$ is the word $\rho(u)$ obtained from $u$ by iteratively deleting factors of the form $a\bar a$ ($a\in \tilde A$). We can then think of $F(A)$ as the set of reduced words on $\tilde A$: the product in $F(A)$ is given by $u \cdot v = \rho(uv)$, and the inverse of $u$ is $\bar u$.

In the sequel, we fix a finite set $A$, with cardinality $r > 1$. If $n\in \N$, we denote by $[n]$ the set of positive integers less than or equal to $n$, and by $\red_n$ (resp. $\red_{\le n}$) the set of reduced words of length $n$ (resp. at most $n$). A reduced word $u$ is called \emph{cyclically reduced} if $u^2$ is reduced, and we let $\cred_n$ (resp. $\cred_{\le n}$) be the set of cyclically reduced words of length $n$ (resp. at most $n$). If $u$ is a reduced word, there exist uniquely defined words $v,w$ such that $w$ is cyclically reduced and $u = v\inv wv$. Then $w$ is called the \emph{cyclic reduction} of $u$, written $\kappa(u)$.

It is easily verified that
\begin{align*}
\red_n = 2r(2r-1)^{n-1} \quad\textrm{and}\quad &2r(2r-1)^{n-2}(2r-2) \le \cred_n \le 2r(2r-1)^{n-1} \\
\red_{\le n} = \Theta((2r-1)^n) \quad\textrm{and}\quad &\cred_{\le n} = \Theta((2r-1)^n).
\end{align*}

If $\vec h = (h_1,\cdots, h_k)$ is a tuple of elements of $F(A)$, we denote by $\langle\vec h\rangle$ the \emph{subgroup of $F(A)$ generated by $\vec h$}: it is the set of all products of the elements of $\vec h$ and their inverses. And we denote by $\llangle\vec h\rrangle$ the normal closure of $\langle\vec h\rangle$, namely the set of all products of conjugates $h_i^g = g\inv h_ig$ of the elements of $\vec h$ ($1\le i\le k$, $g\in F(A)$) and their inverses. The \emph{group presented by the relators $\vec h$}, written $\langle A \mid \vec h\rangle$, is the quotient $F(A)/\llangle\vec h\rrangle$.

If $\vec h = (h_1,\ldots, h_k)$ and $\kappa(\vec h) = (\kappa(h_1), \ldots, \kappa(h_k))$, then $\vec h$ and $\kappa(\vec h)$ present the same group: that is, $\langle A \mid \vec h \rangle = \langle A \mid \kappa(\vec h)\rangle$. It is therefore customary, when considering finite group presentations, to assume that the relators are all cyclically reduced.

In general, if there exists a surjective morphism $\mu\colon F(A) \to G$, we say that $G$ is \emph{$A$-generated}. Then a word $u \in F(A)$ is called \emph{geodesic} if it has minimum length in $\mu\inv(\mu(u))$.

Properties of interest for subgroups of $G$ are, for instance, whether they are free or quasi-convex. Recall that a subgroup $H$ of $G$ is \emph{quasi-convex} if there exists a constant $k > 0$ such that, for every geodesic word $u = a_1\cdots a_n$ such that $\mu(u) \in H$, and for every $1\le i \le n$ there exists a word $v_i$ of length at most $k$ such that $\mu(a_1\cdots a_iv_i) \in H$. 

We are also interested in malnormality and purity: a subgroup $H$ is \emph{almost malnormal} (resp. \emph{malnormal}) if $H^g \cap H$ is finite (resp. trivial) for every $g\not\in H$. Moreover, $H$ is \emph{almost pure} (resp. \emph{pure}, also known as \emph{isolated} or \emph{closed under radical}) if $x^n\in H$ implies $x\in H$ for any $n\ne 0$ and any element $x\in G$ of infinite order (resp.  any $x\in G$). Note that malnormality and almost malnormality (resp. purity and almost purity) are equivalent in torsion-free groups. It is easily verified that an almost malnormal (resp. malnormal) subgroup is almost pure (resp. pure).

It is a classical result that every finitely generated subgroup of a free group is free (Nielsen \cite{1918:Nielsen}) and quasi-convex (Gromov \cite{1987:Gromov}). In addition, it is decidable whether a finitely generated subgroup of a free group is malnormal \cite{1999:BaumslagMyasnikovRemeslennikov} and whether it is pure \cite{2000:BirgetMargolisMeakin}, see Section~\ref{sec: stallings}. In contrast, these properties are not decidable in a general finitely presented group, even if hyperbolic \cite{2001:BridsonWise}. Quasi-convexity is also not decidable in general, even in hyperbolic or small cancellation groups \cite{1982:Rips}. Almost malnormality is however decidable for quasi-convex subgroups of hyperbolic groups \cite[Corollary 6.8]{2016:KharlampovichMiasnikovWeil}.

Finally, let us mention the property of Whitehead minimality for finitely generated subgroups of free groups: we say that $H$ is \emph{Whitehead minimal} if it has minimum size in its automorphic orbit, where the size of a subgroup is defined in terms of its Stallings graph, see Section~\ref{sec: stallings} below. In the case of a cyclic subgroup $H = \langle u\rangle$, if $u = v\inv \kappa(u) v$, then the size of $H$ is $|v|+|\kappa(u)|$. Whitehead minimality plays an important role in the solution of the automorphic orbit problem, to decide whether two subgroups are in the same orbit under the automorphism group of $F(A)$, see \cite{1984:Gersten,1992:Kalajdzievski}.

For finite group presentations, the emphasis can be on combinatorial properties of the presentation, such as small cancellation properties, or on the geometric properties of the given presented group, typically hyperbolicity. One of the main small cancellation properties is Property $C'(\lambda)$ (for some $0 < \lambda < 1$), which is defined as follows. A \emph{piece} in a tuple $\vec h$ of cyclically reduced words is a word $u$ which has at least two occurrences as a prefix of a cyclic conjugate of a word in $\vec h$. A finite presentation $\langle A \mid \vec h\rangle$ \emph{satisfies the small cancellation property $C'(\lambda)$} if a piece $u$ in $\vec h = (h_1,\ldots, h_k)$ satisfies $|u| < \lambda|h_i|$ for every $i$ such that $u$ is a prefix of a cyclic conjugate of $h_i$. This is an important property since it is well known that if $\vec h$ has Property $C'(\frac16)$, then the group $\langle A \mid \vec h \rangle$ is hyperbolic \cite[0.2.A]{1987:Gromov}. Other small cancellation properties are discussed in Section~\ref{sec: density model}.

%%%%%%%%%%%%%%%%%%%
\subsection{Models of randomness}\label{sec: models of randomness}

In this paper, the general model of randomness on a set $S$ which we will consider, consists in the choice of a sequence $(\P_n)_n$ of probability laws on $S$. For instance, the set $S$ could be the set of all $k$-relator presentations (for a fixed value of $k$), that is, the set of all $k$-tuples of cyclically reduced words, and the law $\P_n$ could be the uniform probability law with support the presentations where every relator has length at most $n$.

This general approach covers the classical models considered in the literature, such as the Arzhantseva-Ol'shanski\u\i\ model \cite{1996:ArzhantsevaOlshanskii} or Gromov's density model \cite{1993:Gromov}. It also allows us to consider probability laws that do not give equal weight to words of equal length, see Section~\ref{sec: non-uniform distributions} below.

A subset $X$ of $S$ is \emph{negligible} if the probability for an element of $S$ to be in $X$, tends to 0 when $n$ tends to infinity; that is, if $\lim_n \P_n(X) = 0$. If this sequence converges exponentially fast (that is: $\P_n(X)$ is $\O(e^{-cn})$ for some $c > 0$), we say that $X$ is \emph{exponentially negligible}. The set $X$ is \emph{generic} (resp. \emph{exponentially generic}) if its complement is negligible (resp. exponentially negligible).

%%%%%%%%%%%%%%%%%%%
\section{Random finite presentations}\label{sec: presentations}

\subsection{The density model}\label{sec: density model}
The density model was introduced by Gromov \cite{1993:Gromov}. Let $0 < d < 1$ be a real number. In the density $d$ model, the set $S$ (with reference to the notation in Section~\ref{sec: models of randomness}) is the set of all finite tuples of cyclically reduced words and $\P_n$ is the uniform probability law with support the set of $|\cred_n|^d$-tuples of elements of $\cred_n$. We say that a property is generic (resp. negligible) \emph{at density $d$} if it is generic (resp. negligible) in the density $d$ model.

In this model, small cancellation properties are generic at low enough density. For Property $C'(\lambda)$ ($0 < \lambda < 1$), we have an interesting, so-called \emph{phase transition} statement (\cite[9.B]{1993:Gromov}, see also \cite[Section I.2.a]{2005:Ollivier}).

\begin{theorem}\label{thm: transition for small cancellation}
Let $0 < d < 1$ and $0 < \lambda < \frac12$. If $d < \frac\lambda2$, then at density $d$, a random finite presentation exponentially generically satisfies Property $C'(\lambda)$. If instead $d > \frac\lambda2$, then at density $d$, a random finite presentation exponentially generically does not satisfy $C'(\lambda)$.
\end{theorem}

As noted earlier, if $\vec h$ satisfies Property $C'(\frac16)$, then the group $\langle A \mid \vec h\rangle$ is hyperbolic but the condition is not necessary. Theorem~\ref{thm: transition for small cancellation} shows that, in the density model and up to density $\frac1{12}$, a finitely presented group is exponentially generically hyperbolic. Yet the property holds for higher densities, and we have another phase transition theorem.

Let us say that a finitely presentated group $G = \langle A  \mid \vec h\rangle$, where $\vec h$ consists of cyclically reduced words of equal length $n$, is \emph{degenerate} if $G$ is trivial, or if $G$ is the 2-element group and $n$ is even. Then we have the following result, again a phase transition theorem, due to Ollivier \cite{2005:Ollivier}.

\begin{theorem}\label{thm: transition for hyperbolic}
Let $0 < d < 1$. If $d < \frac12$, then at density $d$, a random finite presentation exponentially generically presents an infinite hyperbolic group. If instead $d > \frac12$, then at density $d$, a random finite presentation exponentially generically presents a degenerate group.
\end{theorem}

The proofs of Theorem~\ref{thm: transition for small cancellation} and of the statement in Theorem~\ref{thm: transition for hyperbolic} about density greater than $\frac12$ reduce to counting arguments on words (see Section~\ref{sec: non-uniform distributions} for a generalization). The proof that hyperbolicity is generic at densities between $\frac1{12}$ and $\frac12$ --- that is: greater than the critical value for Property $C'(\frac16)$ ---, is more complex and involves the combinatorics of van Kampen diagrams.

\begin{remark}\label{rk: ball vs sphere}
A natural variant of the density model considers tuples of words of length at most $n$, instead of words of length exactly $n$. More precisely, $\P_n$ is chosen to be the uniform probability law with support the set of $|\cred_{\le n}|^d$-tuples of words in $\cred_{\le n}$. Ollivier shows in \cite{2005:Ollivier} that Theorems~\ref{thm: transition for small cancellation} and~\ref{thm: transition for hyperbolic} also hold for this model.
\end{remark}

\begin{remark}
The statement on hyperbolicity in Theorem~\ref{thm: transition for hyperbolic} has an important predecessor. For a fixed number $k$ of relators and a fixed $k$-tuple of lengths $(\ell_1,\ldots,\ell_k)$, consider the finite presentations with $k$ relators of length, respectively $\ell_1,\ldots,\ell_k$. The probability that such a presentation presents an infinite hyperbolic group, tends exponentially fast to $1$ when $\min(\ell_i)_{1\leq i \leq k}$ tends to infinity (while $k$ remains fixed). This was originally stated by Gromov \cite{1987:Gromov}, and proved by Champetier \cite{1995:Champetier} and Ol'shanski\u\i\ \cite{1992:Olshanskii}.
\end{remark}

The small cancellation property $C'(\lambda)$ for a tuple of cyclically reduced words $\vec h$ can be interpreted geometrically as follows: in any reduced van Kampen diagram (w.r.t. the presentation $\langle A \mid \vec h\rangle$), a segment of consecutive edges in the boundary between two adjacent faces $f$ and $f'$ (namely, in $\partial f \cap \partial f'$) has length at most $\lambda\min(|\partial f|, |\partial f'|)$. \emph{Greendlinger's property} (as interpreted by Ollivier \cite{2007:Ollivier}) is of the same nature: it states that in any reduced van Kampen diagram $D$ with more than one face, there exist two faces $f$ and $f'$ for which there are segments of consecutive edges of $\partial f \cap \partial D$ (resp. $\partial f' \cap \partial D$) of length at least $\frac12|\partial f|$ (resp. $\frac12|\partial f'|$).

A closely related property of a tuple of relators $\vec h$ is whether Dehn's algorithm works for the corresponding presentation. More precisely, Dehn's algorithm is the following (non-deterministic) process applied to a reduced word $w$: if $w$ if of the form $w = w_1uw_2$ for some word $u$ such that $uv$ is a cyclic permutation of a relator and $|v| < |u|$, then replace $w$ by the reduction of $w_1v\inv w_2$ (which is a shorter word), and repeat. It is clear that this process always terminates, and that if it terminates with the empty word, then $w$ is equal to 1 in the group $G = \langle A \mid \vec h\rangle$. The converse does not hold in general, but we say that $\vec h$ is a \emph{Dehn presentation} if it does, that is, if every reduced word $w$ that is trivial in $G$, contains a factor which is a prefix of some cyclic conjugate $h'$ of a relator, of length greater than $\frac12|h'|$. It is clear that the word problem (given a word $u$, is it equal to 1 in $G$) admits an efficient decision algorithm in groups given by a Dehn presentation. Note that a tuple $\vec h$ with the Greendlinger property provides a Dehn presentation. Moreover, every hyperbolic group has a computable Dehn presentation \cite[Theorem 2.12]{1991:AlonsoBradyCooper} and \cite{1989:Lysenok}.

Greendlinger \cite{1960:Greendlinger} shows that a tuple $\vec h$ with property $C'(\frac16)$ also has (a stronger version of) the Greendlinger property defined above and yields a Dehn presentation. Theorem~\ref{thm: transition for small cancellation} shows that this situation is exponentially generic at density $d < \frac1{12}$. Ollivier proved a phase transition result regarding this property, with critical density higher than $\frac1{12}$ \cite{2007:Ollivier}.

\begin{theorem}\label{thm: transition for Greendlinger}
Let $0 < d < 1$. If $d < \frac15$, then at density $d$, a random finite presentation generically is Dehn and has the Greendlinger property. If instead $d > \frac15$, then at density $d$, a random finite presentation generically fails both properties.
\end{theorem}

\medskip

Ollivier \cite{2004:Ollivier} also considered finite presentations based on a given, fixed hyperbolic $A$-generated group $G$, that is, quotients of $G$ by the normal subgroup generated by a tuple $\vec h$ of elements of $G$, that can be taken randomly. There are actually two ways of generating $\vec h$. Let $\pi\colon F(A) \rightarrow G$ be the canonical onto morphism: then one can draw uniformly at random a tuple $\vec h$ of cyclically reduced words (that is, of elements of $F(A)$) and consider the quotient $G /\!\llangle\pi(\vec h)\rrangle$, or one can draw uniformly at random a tuple of cyclically reduced words that are geodesic for $G$. Here it is useful to remember that a hyperbolic group is geodesically automatic \cite{1992:EpsteinCannonHolt} --- in particular its language of geodesics $L$ is regular ---, and there is a linear time algorithm to randomly generate elements of $L$ of a given length \cite{2012:BernardiGimenez}.

To state Ollivier's result, let us recall the definition of the co-growth of $G$, relative to the morphism $\pi$, under the hypothesis that $\pi$ is not an isomorphism, that is, $G$ is not free over $A$. Let $r = |A|$. Then $\cogrowth(G) = \lim \frac1n \log_{2r-1}(|H_n|)$, where $H_n$ is the set of reduced words of length $n$ in $\ker\pi$ (and the limit is taken over all even values of $n$, to account for the situation where no odd length reduced word is in $\ker\pi$). This invariant of $G$ (and $\pi$) was introduced by Grigorchuk \cite{1980:Grigorchuk}, who proved that it is always greater than $\frac12$ and less than or equal to 1 and, using a result of Kesten \cite{1959:Kesten}, that it is equal to 1 if and only if $G$ is amenable (amenability is an important property which, in the case of discrete groups, is equivalent to the existence of a left-invariant, finitely additive probability measure on $G$). The above definition does not apply if $G$ is free over $A$, but it is convenient to let $\cogrowth(F(A)) = \frac12$ (see \cite[Section 1.2]{2004:Ollivier} for a discussion). In particular, the following elegant phase transition statement generalizes Theorem~\ref{thm: transition for hyperbolic}.

\begin{theorem}\label{thm: quotient of hyperbolic}
Let $G$ be a hyperbolic and torsion-free $A$-generated group and let $\pi\colon F(A)\rightarrow G$ be the canonical mapping. Let $0 < d < 1$. If $d < 1-\cogrowth(G)$, then at density $d$, a random quotient $G /\!\llangle\pi(\vec h)\rrangle$ is exponentially generically hyperbolic. If $d > 1-\cogrowth(G)$, then $G /\!\llangle\pi(\vec h)\rrangle$ is exponentially generically degenerate.

If instead, we take a tuple of cyclically reduced words of length $n$ that are geodesic for $G$, then the phase transition between hyperbolicity and degeneracy is at density $\frac12$. 
\end{theorem}

\begin{remark}
Theorem~\ref{thm: quotient of hyperbolic} above is for torsion-free hyperbolic groups $G$. It actually holds as well if $G$ is hyperbolic and has \emph{harmless torsion}, that is, if every torsion element either sits in the virtual center of $G$, or has a finite or virtually $\Z$ centralizer, see \cite{2004:Ollivier}.
\end{remark}

Finally we note another phase transition theorem, due to \.Zuk, about Kazhdan's property (T) --- a property of the unitary representations of a group --- for discrete groups \cite{2003:Zuk}.  

\begin{theorem}\label{thm: transition for T property}
Let $0 < d < \frac12$. If $d < \frac13$, then at density $d$, a random finite presentation generically does not present a group with Kazhdan's property (T). If instead $d > \frac13$, then at density $d$, a random finite presentation generically presents a group satisfying this property.
\end{theorem}

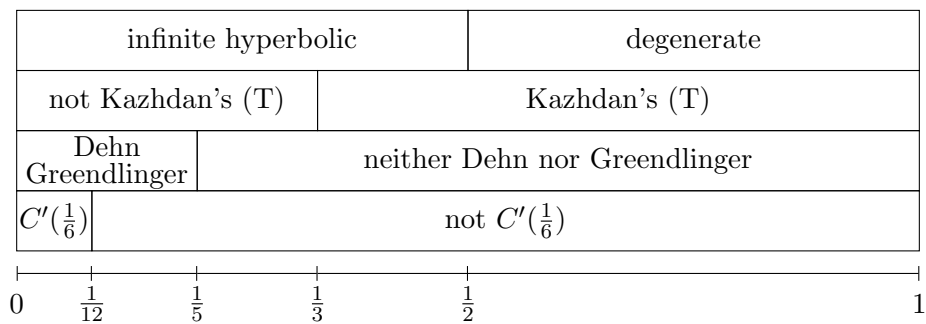
\begin{figure}[htb]
\begin{center}
{\small
\begin{tikzpicture}
\draw (0,0) rectangle (1,.8); \node at (.5,0.4) {$C'(\frac16)$};
\draw (1,0) rectangle (12,.8);\node at (6.5,0.4) {not $C'(\frac16)$};

\draw (0,.8) rectangle (2.4,1.6); \node at (1.2,1.4) {Dehn};\node at (1.2,1) { Greendlinger};
\draw (2.4,.8) rectangle (12,1.6);\node at (7.2,1.2) {neither Dehn nor Greendlinger};

\draw (0,1.6) rectangle (4,2.4); \node at (2,2) {not Kazhdan's (T)};
\draw (4,1.6) rectangle (12,2.4);\node at (8,2) {Kazhdan's (T)};

\draw (0,2.4) rectangle (6,3.2); \node at (3,2.8) {infinite hyperbolic};
\draw (6,2.4) rectangle (12,3.2);\node at (9,2.8) {degenerate};

\node at (0,-0.7) {\small $0$};
\draw[|-|] (0,-.3) -- (1,-.3); \node at (1,-0.7) {\small $\frac1{12}$};
\draw[-|] (1,-.3) -- (2.4,-.3);\node at (2.4,-0.7) {\small $\frac1{5}$};
\draw[-|] (2.4,-.3) -- (4,-.3);\node at (4,-0.7) {\small $\frac1{3}$};
\draw[-|] (4,-.3) -- (6,-.3);\node at (6,-0.7) {\small $\frac1{2}$};
\draw[-|] (6,-.3) -- (12,-.3);\node at (12,-0.7) {\small $1$};
\end{tikzpicture}
}
\end{center}
\caption{\small Phase transitions for properties of random presentations in the density model.}
\end{figure}

%%%%%%%%%%%%%
\subsection{The few relators model}\label{sec: few relators model}

The few relators model was introduced by Arzhantseva and Ol’shanski\u\i\ \cite{1996:ArzhantsevaOlshanskii}.
In this model, the number of relators is fixed, say $k \ge 1$. Then the set $S$ on which we define a model of randomness (see Section~\ref{sec: models of randomness}) is the set $\cred^k$ of all $k$-tuples of cyclically reduced words in $F(A)$ and $\P_n$ is the uniform probability law with support $(\cred_{\le n})^k$.

Observe that if a tuple $\vec h$ of cyclically reduced words satisfies the small cancellation property $C'(\lambda)$ and if $\vec g$ is a sub-tuple of $\vec h$ (that is, the words in $\vec g$ are also in $\vec h$), then $\vec g$ satisfies Property $C'(\lambda)$ as well. From this observation and Theorem~\ref{thm: transition for small cancellation} (actually its variant in Remark~\ref{rk: ball vs sphere}) we deduce the following result, due to Arzhantseva and Ol'shanskii \cite{1996:ArzhantsevaOlshanskii} (see also \cite[Theorem B]{2005:KapovichSchupp}).

\begin{corollary}\label{cor: hyperbolic in few relators}
In the few relator model, a random tuple exponentially generically satisfies Property $C'(\frac16)$ and presents an infinite hyperbolic group.
\end{corollary}

Arzhantseva and Ol’shanski\u\i\ showed further that, in the few relator model, the finitely generated subgroups of a random $k$-relator subgroup are usually free or have finite index \cite{1996:ArzhantsevaOlshanskii}, statements (1) and (2) of Theorem~\ref{thm: subgroups are free} below. Statement (3) is due to Kapovich and Schupp \cite[Theorem B]{2005:KapovichSchupp}. Recall that a \emph{Nielsen move} on a tuple $(x_1,\ldots, x_k)$ of elements of a group $G$ consists in replacing $x_i$ by $x_i\inv$, exchanging $x_i$ and $x_j$ or replacing $x_i$ by $x_ix_j$ for some $i\ne j$. We say that two $k$-tuples are \emph{Nielsen-equivalent} if one can go from one to the other by a sequence of Nielsen moves.

\begin{theorem}\label{thm: subgroups are free}
Let $k,\ell \ge 1$ be integers. In the few relators model with $k$ relators, exponentially generically,
\begin{itemize}
\item[(1)] every $\ell$-generated subgroup of an $A$-generated group $G$, has finite index or is free;

\item[(2)] if $\ell < |A|$, every $\ell$-generated subgroup of $G$ is free and quasi-convex in $G$;

\item[(3)] an $|A|$-tuple which generates a non-free subgroup of $G$ is Nielsen-equivalent to $A$ in $G$. In particular, an $|A|$-tuple which generates a non-free subgroup generates $G$ itself, and every automorphism of $G$ is induced by an automorphism of $F(A)$.
\end{itemize}
\end{theorem}

Arzhantseva also established the following related result  \cite{2000:Arzhantseva}, which refines in a sense Theorem~\ref{thm: subgroups are free}~(1). Here a set of generators for the $\ell$-generated subgroup of $G$ is fixed in advance (as a tuple of words in $F(A)$), and it is assumed that it generates an infinite-index subgroup of $F(A)$.

\begin{theorem}
Let $H$ be a finitely generated, infinite index subgroup of $F(A)$. In the few relators model with $k$ relators, exponentially generically, a finite presentation $G = \langle A \mid \vec h\rangle$ ($\vec h \in (\cred_{\le n})^k$) is such that the canonical morphism
$\phi\colon F(A) \rightarrow G$ is injective on $H$ (so $\phi(H)$ is free) and $\phi(H)$ has infinite index in $G$.
\end{theorem}

We also note that Kapovich and Schupp extended Theorem~\ref{thm: subgroups are free} to the density model \cite{2008:KapovichSchupp}, with density bounds that depend on both parameters $k$ and $\ell$.

\begin{theorem}
Let $A$ be a fixed alphabet. For every $k,\ell \ge 1$, there exists $0 < d(k,\ell) < 1$ such that, at every density $d < d(k,\ell)$, generically, an $\ell$-generated subgroup of an $A$-generated group presented by a random $k$-tuple of relators has finite index or is free.

Also, for every $k \ge 1$, there exists $0 < d(k) < 1$ such that, at every density $d < d(k)$, every $(k - 1)$-generated subgroup of an $A$-generated group presented by a random $k$-tuple of relators is free. But there is no single value of $d$ such that this holds independently of $k$ (that is: $\lim_{k\to\infty} d(k) = 0$).
\end{theorem}

%%%%%%%%%%%%%%%%%%
\subsection{1-relator groups}\label{sec: 1-relator}

If $u$ is a cyclically reduced word, let $G_u = \langle A \mid u\rangle$. 1-relator groups are of course covered by the few relators model, and the results of Section~\ref{sec: few relators model} apply to them. But more specific results are known for random 1-relator presentations.

Magnus showed that if $u,v\in F(A)$, then the normal closures of the subgroups generated by $u$ and $v$, written $\llangle u\rrangle$ and $\llangle v\rrangle$, are equal if and only if $u$ is a conjugate of $v$ or $v\inv$ (see \cite[Prop. II.5.8]{1977:LyndonSchupp}). Kapovich and Schupp combine this with Theorem~\ref{thm: subgroups are free}~(3) to show the following \cite[Theorem A]{2005:KapovichSchupp}.

\begin{theorem}\label{thm: Magnus and KS}
There exists an exponentially generic (and decidable) class $P$ of cyclically reduced words such that, if $u,v\in P$, then $G_u$ and $G_v$ are isomorphic if and only if there exists an automorphism $\phi$ of $F(A)$ such that $\phi(u) \in \{v,v\inv\}$. In particular, the isomorphism problem for 1-relator groups with presentation in $P$ is decidable.
\end{theorem}

We now explain how this result gives access to generic properties of (isomorphism classes of) 1-relator groups and not just of 1-relator presentations. This is a more explicit rendering of arguments which can be found in particular in Ollivier \cite[Section II.3]{2005:Ollivier}, Kapovich, Schupp and Shpilrain \cite{2006:KapovichSchuppShpilrain} and Sapir and \v{S}pakulová \cite{2011:SapirSpakulova}. For this discussion we consider probability laws $\P_n$ for 1-relator presentations and probability laws $\Q_n$ for isomorphism classes of 1-relator groups. More specifically, $\P_n$ is the uniform probability law with support $\cred_{\le n}$ (that is: the probability law for the few relator model, with $k=1$ relator); and $\Q_n$ is the uniform probability law with support the set $T_n$ of isomorphism classes of groups $G_u$ with $|u|\le n$. We let $T = \bigcup_{n\ge 1} T_n$, that is, $T$ is the set of isomorphism classes of 1-relator groups.

Let $H$ be the group of length-preserving automorphisms of $F(A)$, that is, the automorphisms which permute $\tilde A$. Note that $|H| = 2^rr!$, where $r = |A|$. Let also $W$ be the set of strictly Whitehead minimal words, that is, cyclically reduced words $u$ such that $|\phi(u)| > |u|$ for every automorphism $\phi\in \Aut(F(A)) \setminus H$. Kapovich \emph{et al.} \cite{2006:KapovichSchuppShpilrain} show that $W$ is exponentially generic (see \cite{2016:BassinoNicaudWeil} and Section~\ref{sec: Whitehead minimality} for a more general result).

Fix an arbitrary order on $\tilde A$. For each word $u\in\cred$, we let $\tau(u)$ be the lexicographically least element of the set of all cyclic permutations of images of $u$ and $u\inv$ by an automorphism in $H$. If $P$ is the set in Theorem~\ref{thm: Magnus and KS}, we observe that $u\in P$ if and only if $\tau(u) \in P$ (see for instance the description of $P$ in \cite[Section 4]{2005:KapovichSchupp}). The same clearly holds for $W$, and we have $2r \le |\tau\inv(\tau(u))| \le 2^{r+1}|u|r!$ --- where the lower bound corresponds to a word of the form $u = a^{|u|}$. It is immediate that $G_u = G_{\tau(u)}$. Moreover, in view of Theorem~\ref{thm: Magnus and KS}, if $u,v\in P\cap W$, then $G_u$ and $G_v$ are isomorphic if and only if $\tau(u) = \tau(v)$.

\begin{proposition}\label{prop: statistics for 1-relator groups}
Let $X$ be a property of isomorphism classes of 1-relator groups, that is, $X$ is a subset of $T$. Let $Y = \{u\in \cred \mid G_u \in X\}$. If $\P_n(Y) = o(n\inv)$ (resp. $Y$ is exponentially negligible), then $X$ is negligible (resp. exponentially negligible). The same statement holds for genericity instead of negligibility.
\end{proposition}

\proof
Let $Z$ be the set of 1-relator groups $G_u$ such that $u\in W\cap P$, where $P$ is the set in Theorem~\ref{thm: Magnus and KS} and $W$ is the set of strictly Whitehead minimal cyclically reduced words. Since $W \cap P$ is exponentially generic in $\cred$, there exist constants $C,c> 0$ such that $\P_n(W\cap P) \ge 1-Ce^{-cn}$. We have
$$\Q_n(X) = \Q_n(X \cap Z) + \Q_n(X \setminus Z) \le \Q_n(X \cap Z) + \Q_n(T \setminus Z).$$
We first deal with $\Q_n(T \setminus Z)$. Let $\alpha_n = |T_n \cap Z|$ and $\beta_n = |T_n \setminus Z|$. Then $\Q_n(T \setminus Z) = \frac{\beta_n}{\alpha_n+\beta_n}$. Note that $T_n \setminus Z \subseteq \{G_u \mid u \in \cred_{\le n}\setminus(W\cap P)\}$. So $\beta_n \le |\cred_{\le n}\setminus(W\cap P)| \le Ce^{-cn}|\cred_{\le n}|$.

On the other hand, $T_n \cap Z$ is in bijection with $\{\tau(u) \mid u\in \cred_{\le n} \cap (W\cap P)\}$, and it follows that
$$\alpha_n \ge \frac1{2^{r+1}nr!} |\cred_{\le n}\cap (W\cap P)| \ge \frac{1-Ce^{-cn}}{2^{r+1}nr!} |\cred_{\le n}|.$$
Therefore
$$\Q_n(T \setminus Z) = \frac{\beta_n}{\alpha_n+\beta_n} \le \frac{\beta_n}{\alpha_n} \le \frac{Ce^{-cn}}{1-Ce^{-cn}} 2^{r+1}nr!,$$
which vanishes exponentially fast.

Let us now consider $\Q_n(X \cap Z)$. We have
\begin{align*}
\Q_n(X \cap Z) &= \frac{|X \cap Z \cap T_n|}{\alpha_n + \beta_n}\\
&\le \frac{|\{G_u \mid u\in \cred_{\le n} \cap W\cap P, u\in Y\}|}{\alpha_n}\\
&\le \frac{2^{r+1}nr!\ \P_n(Y\cap W \cap P)\ |\cred_{\le n}|}{2r(1-Ce^{-cn})\ |\cred_{\le n}|}\\
&\le 2^r(r-1)!\ \frac{n\P_n(Y)}{1-Ce^{-cn}},
\end{align*}
and this concludes the proof.
\eop
Then the results of Section~\ref{sec: few relators model} (Corollary~\ref{cor: hyperbolic in few relators}, Theorem~\ref{thm: subgroups are free}), together with Proposition~\ref{prop: statistics for 1-relator groups} yield the following.

\begin{corollary}
Exponentially generically, a 1-relator group $G$ is infinite hyperbolic, every automorphism of $G$ is induced by an automorphism of $F(A)$, and every $\ell$-generated subgroup is free and quasi-convex if $\ell < |A|$.
\end{corollary}

Kapovich, Schupp and Shpilrain use the ideas behind Proposition~\ref{prop: statistics for 1-relator groups} to compute an asymptotic equivalent of the number of (isomorphism classes of) 1-relator groups in $T_n$  \cite{2006:KapovichSchuppShpilrain}.

\begin{theorem}
Let $I_n(A)$ be the number of isomorphism classes of 1-relator groups of the form $\langle A \mid u\rangle$ with $|u| \le n$. If $|A| = r$, then $I_n(A)$ is asymptotically equivalent to $\frac1{2^{r+1}r!}\ \frac{(2r-1)^n}{n}$.
\end{theorem}

Finally we note the following result of Sapir and \v Spakulová \cite{2011:SapirSpakulova}. Recall that a group $G$ is \emph{residually $\calP$} (for some property $\calP$) if for all distinct elements $x, y \in G$, there exists a morphism $\phi$ from $G$ to a group having property $\calP$, such that $\phi(x) \ne \phi(y)$. We let \emph{finite-$p$} be the property of being a finite $p$-group. A group $G$ is \emph{coherent} if every finitely generated subgroup is finitely presented.

\begin{theorem}
Suppose that $|A| \ge 3$. Then an $A$-generated 1-relator group is generically residually finite, residually finite-$p$ and coherent.
\end{theorem}

%%%%%%%%%%%%%%%%%%
\subsection{Rigidity properties}\label{sec: rigidity}

Theorem~\ref{thm: Magnus and KS} above gives a generic rigidity property: at least on a large (exponentially generic) set of words, the isomorphism class of the 1-relator group $G_u = \langle A \mid u\rangle$ is uniquely determined by $u$, up to inversion and an automorphism of $F(A)$. That is, the only words $v$ such that $G_v$ is isomorphic to $G_u$ are those that come immediately to mind. As indicated, this result follows from a theorem of Magnus which states that the normal closure $\llangle u\rrangle$ has essentially only one generator as a normal subgroup: $\llangle u\rrangle = \llangle v\rrangle$ if and only if $u$ is a conjugate of $v$ or $v\inv$.

There is no such general statement for normal subgroups generated by a $k$-tuple with $k\ge 2$. A closely related result due to Greendlinger generalizes Magnus's statement, but only for tuples that satisfy the small cancellation property $C'(\frac16)$ \cite{1961:Greendlinger}: if $\vec g$ and $\vec h$ are such tuples, respectively a $k$-tuple and an $\ell$-tuple, and if $\llangle \vec g \rrangle = \llangle \vec h \rrangle$, then $k = \ell$ and there is a re-ordering $\vec g'$ of $\vec g$ such that, for each $i$, $h_i$ is a cyclic permutation of $g'_i$ or ${g'_i}\inv$. The restriction to tuples satisfying $C'(\frac16)$ prevents us from proceeding as in Section~\ref{sec: 1-relator} to prove a more general analogue of Theorem~\ref{thm: Magnus and KS}. Whether Theorem~\ref{thm: Magnus and KS} can be extended to $m$-tuples of cyclically reduced words, is essentially the \emph{Stability Conjecture} formulated by Kapovich and Schupp \cite[Conjecture 1.2]{2005:KapovichSchuppHelv}.

Nevertheless, Kapovich and Schupp show that one can circumvent this obstacle when considering the quotients of the \emph{modular} group $M = \PSL(2,\Z) = \langle a,b \mid a^2, b^3\rangle$. If $\vec h$ is a tuple of cyclically reduced words in $F(a,b)$, we denote by $M_{\vec h}$ the quotient of $M$ by the images of the elements of $\vec h$ in $M$, that is, $M_{\vec h} = \langle a,b \mid a^2, b^3, \vec h\rangle$. Let $\eta$ be the automorphism of $M$ which fixes $a$ and maps $b$ to $b\inv = b^2$. Then the following holds \cite[Theorem A and Corollary 2.5]{2009:KapovichSchupp}.

\begin{theorem}\label{thm: rigidity modular}
For each $k\ge 1$, there exists an exponentially generic (in the $k$-relator model), decidable subset $Q_k$ of $\cred^k$ such that the following holds.
\begin{itemize}
\item If $\vec h\in Q_k$, then the group $M_{\vec h}$ is hyperbolic and one-ended, the generators $a$ and $b$ have order 2 and 3, respectively in $M_{\vec h}$, and all the automorphisms of $M_{\vec h}$ are inner.

\item If $\vec g, \vec h\in Q_k$ and $M_{\vec g}$ and $M_{\vec h}$ are isomorphic, then there is a re-ordering $\vec g'$ of $\vec g$ and a value $\epsilon\in\{0,1\}$ such that, for each $1\le i \le k$, $h_i$ is a cyclic permutation of $\eta^\epsilon(g'_i)$ or $\eta^\epsilon({g'_i}\inv)$.

\item If $\vec g \in Q_k$, $\vec h\in Q_\ell$ are such that the $g_i$ and the $h_j$ all have the same length, and if $M_{\vec g}$ and $M_{\vec h}$ are isomorphic, then $k = \ell$.
\end{itemize}
In the $k$-relator model, the isomorphism problem for quotients of $M$ is exponentially generically solvable in time $\O(n^4)$.
\end{theorem}

The last statement of this theorem is all the more interesting as the isomorphism problem, and even the triviality problem, for quotients of $M$ is undecidable in general (Schupp \cite{1976:Schupp}).

As in Section~\ref{sec: 1-relator}, Theorem~\ref{thm: rigidity modular} can be used to discuss asymptotic properties of $k$-relator quotients of the modular group, rather than of $k$-tuples of relators. The few-relator model for the quotients of $M$ considers the set $T$ of isomorphism classes of $k$-relator quotients of $M$, and the probability laws $\Q_n$ which are uniform on the set $T_n$ of isomorphism classes of groups $M_{\vec h}$ with $\vec h \in (\cred_{\le n})^k$. We can reason as for Proposition~\ref{prop: statistics for 1-relator groups}, modifying the map $\tau$ in such a way that $\tau(h)$ is the lexicographically least element of $h$, $h\inv$, $\eta(h)$ and $\eta(h\inv)$. Then, with essentially the same proof as Proposition~\ref{prop: statistics for 1-relator groups}, we get the following result.

\begin{proposition}\label{prop: statistics for quotients of M}
Let $k\ge 1$ and let $X$ be a property of isomorphism classes of k-relator quotients of the modular group, that is, $X$ is a subset of $T$. Let $Y = \{\vec h\in (\cred)^k \mid M_{\vec h} \in X\}$. If $\P_n(Y) = o(n^{-k})$ (resp. $Y$ is exponentially negligible), then $X$ is negligible (resp. exponentially negligible). The same statement holds for genericity instead of negligibility.
\end{proposition}

As in Section~\ref{sec: 1-relator} again, one can derive from Theorem~\ref{thm: rigidity modular} an asymptotic equivalent of the number of isomorphism classes of $k$-relator quotients of the modular group \cite[Theorem C]{2009:KapovichSchupp}.

\begin{corollary}
Let $k\ge 1$. The number of isomorphism classes of quotients of $M$ by $k$ relators which are cyclically reduced words of length $n$, is asymptotically equivalent to
$$\frac{(2^{\frac n2+1})^k}{2 k! (2n)^k}.$$
\end{corollary}

Kapovich and Schupp go on to give further generic rigidity properties of homomorphisms between quotients of $M$, which are proved to be generically hopfian and co-hopfian (that is, every surjective (resp. injective) endomorphism is an isomorphism), and on the generic incompressibility of the presentations by $k$ relators \cite[Theorems B and D]{2009:KapovichSchupp}.

%%%%%%%%%%%%%%%%%%
\subsection{Nilpotent groups}

We conclude this section with recent results on random groups in a particular class, that of nilpotent groups. If $G$ is a group, the \emph{lower central series} of $G$ is defined by letting $G_1 = G$ and, for $n\ge 1$, $G_{n+1} = [G_n,G]$. That is: $G_{n+1}$ is the subgroup generated by the commutators $[g,h] = g\inv h\inv gh$, with $g\in G_n$ and $h\in G$. Then each $G_n$ is normal in $G$ and $G_{n+1}$ is contained in $G_n$. The group $G$ is said to be \emph{nilpotent of class $s$} if $G_{s+1} = 1$.  In particular, $G_2$ is the derived subgroup of $G$ and the class 1 nilpotent groups are exactly the abelian groups. Nilpotent groups of class 2 are those in which the derived subgroup lies in the center of the group.

Let us extend the commutator notation by letting, for $s\ge 2$, $[x_1, \ldots, x_{s+1}] = [[x_1, \ldots, x_s],x_{s+1}]$. One can show that the class of nilpotent groups of class $s$ is defined by the identity $[x_1, \ldots, x_{s+1}] = 1$. As a result, this class constitutes a variety, and we denote by $\Nil_s(A)$ its free object over the finite alphabet $A$: $\Nil_s(A) = F(A) / F(A)_{s+1}$.

Note that a torsion-free non-cyclic nilpotent group contains a free abelian group of rank 2, a standard obstacle for hyperbolicity: so torsion-free non-cyclic nilpotent groups are not hyperbolic. In particular, they form a negligible set in the few-relator as well as in the density models discussed in the previous sections, and we can not use earlier results to discuss random nilpotent groups. This difficulty was circumvented in several different ways in the literature.

Cordes et al. view finitely presented nilpotent groups as quotients of free nilpotent groups (of a fixed class and rank) by a random tuple of relators whose length tends to infinity \cite{2016:CordesDuchinDuong}. In this model, relators are words over the symmetrized alphabet $\tilde A$. Depending on the number of relators, this extends the few relator and the density models. Garreta et al. extend in \cite{2016:Garreta,2016:GarretaMiasnikovOvchinnikov} the study initiated in \cite{2016:CordesDuchinDuong}. The following result is a summary of \cite[Theorem 29, Proposition 30 and Corollaries 32 and 35]{2016:CordesDuchinDuong} and of \cite[Theorems 3.7 and 4.1]{2016:GarretaMiasnikovOvchinnikov}.

\begin{theorem}\label{thm: Cordes et al}
Let $s\ge 1$, $r \ge 2$, let $A$ be an alphabet of cardinality $r$, let $\Nil_{s,r} = \Nil_s(A)$ be the free nilpotent group of class $s$ over $A$, and let $\pi$ be the canonical morphism from $\tilde A^*$ onto $\Nil_{s,r}$. 

In the density model, at any density $d> 0$, a random quotient $\Nil_{s,r}/\llangle\pi(\vec h)\rrangle$ is generically trivial. In fact, this holds in any model where the size of the tuple of relators is not bounded.

In the few relator model with $k$ relators with $k \le r-2$, a random quotient $\Nil_{s,r}/\llangle\pi(\vec h)\rrangle$ is generically non abelian and regular (that is: every element of the center of $G$ has a non-trivial power in the derived subgroup).

If $k = r-1$, then such a quotient is generically virtually abelian (it has an abelian finite index subgroup), and if $k = r$, then it is generically finite. In either case, it is abelian if and only if it is cyclic. Finally, if $k \ge r+1$, then it is generically finite and abelian.

In the particular case where $r = 2$, $k = 1$ and $s\ge 2$, the probability that a random 1-relator quotient of $\Nil_{s,2}$ is cyclic (and hence, abelian) tends to $\frac6{\pi^2}$.
\end{theorem}

Cordes et al. also give a full classification of the 1-relator quotients of $\Nil_{2,2}$ (the Heisenberg group) \cite[Section 3]{2016:CordesDuchinDuong}. Moreover, they deduce from Theorem~\ref{thm: Cordes et al} the following result on random finitely presented groups \cite[Corollary 36]{2016:CordesDuchinDuong}.

\begin{corollary}
In the density model, at any density $d> 0$, and in any model where the size of the tuple of relators is not bounded, a random tuple $\vec h$ generically presents a perfect group (that is: a group $G$ such that $[G,G] = G$, or equivalently, a group whose abelian quotient is trivial).
\end{corollary}

Delp \emph{et al.} use a different view of nilpotent groups \cite{2016:DelpDymarzSchaffer-Cohen}: it is well known that every torsion-free nilpotent group embeds in $U_n(\Z)$ for some $n \ge 2$, where $U_n(\Z)$ is the group of upper-triangular matrices with entries in $\Z$ and diagonal elements equal to $1$. If $1 \le i < n$, let $a_{i,n}$ be the matrix in $U_n(\Z)$ with coefficients 1 on the diagonal and on row $i$ and column $i+1$, and all other coefficients $0$. Then $A_n = \{a_{1,n},\ldots,a_{n-1,n}\}$ generates $U_n(\Z)$. Let $\ell\colon\N\to\N$ be a function such that $\lim\ell(n) = \infty$ when $n$ tends to infinity. We let $G_{\ell,n}$ be the subgroup of $U_n(\Z)$ generated by a random pair of words of length $\ell(n)$ on alphabet $\tilde A_n$: in the language of Section~\ref{sec: models of randomness}, $S$ is the set of pairs of words on an alphabet of the form $\tilde A_n$ for some $n\ge 2$, $S_n$ is the set of all pairs of length $\ell(n)$ words on alphabet $\tilde A_n$, and $\P_n$ is the uniform probability law with support $S_n$. Then we have the following result, a combination of \cite[Theorems 1 and 2]{2016:DelpDymarzSchaffer-Cohen}. Note that $U_n(\Z)$ is nilpotent of class $n-1$: we say that a subgroup of $U_n(\Z)$ has full class if it is nilpotent of class $n-1$.

\begin{theorem}
Let $\ell\colon\N\to\N$ be a function such that $\lim\ell(n) = \infty$.
\begin{itemize}
\item If $\ell = o(\sqrt n)$, then $G_{\ell,n}$ is generically abelian (that is: of class 1). If $\sqrt n = o(\ell(n))$, then $G_{\ell,n}$ is generically non abelian. And if $\ell(n) = c\sqrt n$, then the probability that $G_{\ell,n}$ is abelian tends to $e^{-2c^2}$.

\item If $\ell = o(n^2)$, then $G_{\ell,n}$ generically does not have full class; and if $n^3 = o(\ell(n))$, then $G_{\ell,n}$ generically has full class.
\end{itemize}
\end{theorem}

Garreta \emph{et al.} use yet another representation of nilpotent groups \cite{2016:Garreta,2016:GarretaMiasnikovOvchinnikov-2}, the polycyclic presentation. A group $G$ is \emph{polycyclic} if it admits a sequence of subgroups $1 = H_n \le H_{n-1} \le \cdots \le H_1 = G$ such that, for every $1 < i \le n$, $H_i$ is normal in $H_{i-1}$ and $H_{i-1}/H_i$ is cyclic. It is elementary to verify that every finitely generated nilpotent group is polycyclic. Polycyclic groups admit presentations of a particular form, the so-called \emph{polycyclic presentations} (see \cite{2005:HoltEickOBrien} for a precise description), which can be characterized by a $k$-tuple of integers, where $k$ is a function of the number of generators in the presentation. 

In the case of torsion-free nilpotent groups, polycyclic presentations with generators $x_1,\ldots, x_r$ have relators of the following form, called a \emph{torsion-free nilpotent presentation}:
\begin{align*}
[x_j,x_i] &= x_{j+1}^{b_{i,j,j+1}}\cdots x_r^{b_{i,j,r}} \\
[x_j,x_i\inv] &= x_{j+1}^{c_{i,j,j+1}}\cdots x_r^{c_{i,j,r}},
\end{align*}
for all $1 \le i < j \le r$, where the $b_{i,j,h}$ and $c_{i,j,h}$ ($1 \le i < j < h\le r$) are integers. Garreta \emph{et al.} introduce a notion of random torsion-free nilpotent presentations as follows \cite{2016:Garreta,2016:GarretaMiasnikovOvchinnikov-2}: with the number $r$ of generators fixed, they let $S$ be the set of tuples $(b_{i,j,h},c_{i,j,h})_{1 \le i < j < h\le r}$ of integers, $S_n$ be the set of those tuples whose components sit in the interval $[-n,n]$ and $\P_n$ be the uniform law with support $S_n$. They then show the following result \cite[Lemma 8]{2016:Garreta,2016:GarretaMiasnikovOvchinnikov-2}.

\begin{proposition}
Let $r\ge 2$. The group presented by a random torsion-free nilpotent presentation is generically finite.
\end{proposition}

The situation becomes more interesting if one restricts one's attention to torsion-free nilpotent groups of class 2, also known as \emph{$\tau_2$-groups}. Recall that these are the torsion-free groups where the derived subgroup is contained in the center. In particular, the derived subgroup and the center are both free abelian groups. In the case of $\tau_2$-groups, torsion-free nilpotent presentations can be simplified to the following \emph{$\tau_2$-presentations}:
$$\langle A, C \mid [a_i,c_h] = 1, [c_h,c_k] = 1, [a_i,a_j] = \prod_{1\le h\le m} c_h^{\alpha_{i,j,h}},\textrm{ $1\le i < j \le \ell$}\rangle,$$
for some $A = \{a_1,\ldots, a_\ell\}$, $C = \{c_1, \ldots, c_m\}$ ($\ell,m\ge 0$) and some choice of $\alpha_{i,j,h}\in \Z$ ($1\le i,j\le \ell$ and $1\le h\le m$). If $\vec\alpha = (\alpha_{i,j,h})_{i,j,h}$, we denote the group thus presented by $G(A,C,\vec\alpha)$. Note that $C$ generates a free abelian group of rank $m$, contained in the center of $G(A,C,\vec\alpha)$.

For a fixed choice of $A$ and $C$, a natural notion of randomness is given by letting $S$ be the set of all the tuples $\vec\alpha$ of the appropriate size, $S_n$ be the set of these tuples where every element has absolute value at most $n$, and $\P_n$ be the uniform probability law with support $S_n$. In this situation, Garreta \emph{et al.} show the following \cite[Theorems 4 and 5]{2016:Garreta,2016:GarretaMiasnikovOvchinnikov-2}.

\begin{theorem}
Let $\ell, m \ge 0$ and let $G$ be the group presented by a random $\tau_2$-presentation on the pair of alphabets $A = \{a_1,\ldots, a_\ell\}$ and $C = \{c_1, \ldots, c_m\}$.

If $\ell-1 \le m$, then generically $C$ generates $Z(G)$, the center of $G$, and $G$ is directly indecomposable into non-abelian factors.

If $m \le \frac{\ell(\ell-1)}2$, then generically the derived subgroup $G'$ of $G$ has finite index in $Z(G)$. If $\ell-1 \le m \le \frac{\ell(\ell-1)}2$, then $G$ is generically regular.

If $m > \frac{\ell(\ell-1)}2$, then $G$ is not regular, $G'$ is freely generated (as an abelian group) by the $[a_i,a_j]$ ($1\le i < j \le \ell$), and $G'$ generically has infinite index in $Z(G)$.
\end{theorem}

Garreta \emph{et al.} also discuss whether the diophantine problem is generically decidable in a random $\tau_2$-group \cite{2016:GarretaMiasnikovOvchinnikov-2}.

%%%%%%%%%%%%%%%%%%%
\section{Random subgroups}\label{sec: random subgroups}

Before we discuss asymptotic properties of finitely generated subgroups, let us introduce a privileged tool to describe and reason about subgroups of free groups. Most of this section is devoted to this type of subgroups, only Section~\ref{sec: subgroups of other groups} below goes beyond the free group case.

%%%%%%%%%%%%%%%%
\subsection{Stallings graph of a subgroup}\label{sec: stallings}

It is classical to represent the finitely generated subgroups of a free group by a finite labeled graph, subject to certain combinatorial constraints. An \emph{$A$-graph} is a finite graph $\Gamma$ whose edges are labeled by elements of $A$. It can be seen also as a transition system on alphabet $\tilde A$, with the convention that every $a$-edge from $p$ to $q$ represents an $a$-transition from $p$ to $q$ and an $\bar a$-transition from $q$ to $p$. Say that $\Gamma$ is \emph{reduced} if
it is connected and if no two edges with the same label start (resp. end) at the same vertex: this is equivalent to stating that the corresponding transition system is deterministic and co-deterministic.  If 1 is a vertex of $\Gamma$, we say that $(\Gamma,1)$ is \emph{rooted} if every vertex, except possibly 1, has valency at least 2.

\begin{figure}[ht]
\begin{minipage}{.24\textwidth}
\centering
\includegraphics[scale=.9]{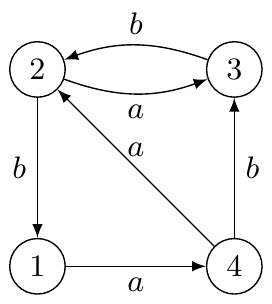}
\end{minipage}
\begin{minipage}{.75\textwidth}
\caption{\small The Stallings graph of $H = \langle aab, ab\bar a b, abbb\rangle$. The reduced word $u=aa\bar b\bar a b$ is
in $H$ as it is accepted by $\Gamma(H)$: it labels a path starting from $1$ and ending at $1$, with edges being
used backward when reading a negative letter. Since every vertex has valency at least 2, this graph is cyclically reduced.
\label{fig:stallings}}
\end{minipage}
\end{figure}

We say that $\Gamma$ is \emph{cyclically reduced} if it is reduced and every vertex has valency at least 2. The $A$-graph in Fig.~\ref{fig:stallings} is cyclically reduced. If $\Gamma$ is reduced, we denote by $\kappa(\Gamma)$ the \emph{cyclic reduction of $\Gamma$}, namely the cyclically reduced $A$-graph obtained from $\Gamma$ by iteratively removing vertices of valency 1 and the edges adjacent to them. 

If $H$ is a finitely generated subgroup of $F(A)$, there exists a unique reduced rooted graph $(\Gamma(H),1)$, called the \emph{Stallings graph} of $H$, such that $H$ is exactly the set of reduced words accepted by $(\Gamma(H),1)$: a reduced word is accepted when it labels a loop starting and ending at $1$. Moreover, this graph can be effectively computed given a tuple of reduced words generating $H$, in time $\O(n\log^*n)$ (that is to say: almost linear) \cite{1983:Stallings,2006:Touikan}.

Conversely, every rooted reduced $A$-graph ($\Gamma,1)$ is the Stallings graph of a (computable) finitely generated subgroup of $F(A)$. Moreover, two subgroups $H_1$ and $H_2$ are conjugated if and only if $\kappa(\Gamma(H_1)) = \kappa(\Gamma(H_2))$. Thus the reduced $A$-graph $\kappa(\Gamma(H))$ can be seen as a representative of the conjugacy class of $H$.

Many interesting properties of a subgroup $H$ can be characterized -- and often decided -- in terms of the Stallings graph $\Gamma(H)$. We list here a couple that will be used in the sequel:
\begin{itemize}
\item  the rank of $H$ is equal to $|E|-|V|+1$, where $E$ (resp. $V$) is the number of edges (resp. vertices) of $\Gamma(H)$ \cite{1983:Stallings};

\item $H$ has finite index in $F(A)$ if and only if $\Gamma(H)$ has the maximal number of possible edges, namely $|V||A|$ (one edge starting from every vertex, labeled by every letter of $A$), and if this is the case, then the index of $H$ in $F(A)$ is $|V|$ \cite{1983:Stallings};

\item $H$ is malnormal if and only if no non-empty reduced word $u$ labels a loop in $\Gamma(H)$ at two different vertices, if and only if every non-diagonal connected component of the direct product $\Gamma(H) \times \Gamma(H)$ (in the category of $A$-graphs) is a tree \cite{1999:BaumslagMyasnikovRemeslennikov};

\item $H$ is pure if and only if $\Gamma(H)$ never has a loop labeled $u^n$ ($u$ a non-empty reduced word, $n\ge 0$) at a vertex $v$, without having in fact a $u$-labeled loop at that vertex \cite{2000:BirgetMargolisMeakin}.
\end{itemize}

%%%%%%%%%%%%%%%%
\subsection{The central tree property and its consequences}\label{sec: ctp}

Let $\vec h = (h_1,\ldots, h_k)$ be a a tuple of reduced words in $F(A)$ and let $\vec h^\pm$ be the $2k$-tuple consisting of the elements of $\vec h$ and their inverses. Let $\min(\vec h) = \min_i|h_i|$ and let $\lcp(\vec h)$ be the length of the longest common prefix of the elements of $\vec h^\pm$.

We say that $\vec h$ has the \emph{central tree property} if $2\lcp\vec h < \min\vec h$. This property is identified explicitly by Bassino \emph{et al.} \cite[Section 1.3]{2016:BassinoNicaudWeilb}, but it is underlying the reasoning in the work of Arzhantseva and Ol'shanski\u\i\ \cite{1996:ArzhantsevaOlshanskii}, Jitsukawa \cite{2002:Jitsukawa} and several others. The central tree property, which could also be termed a small initial cancellation property, has the following interesting consequences. All are easily verified, except perhaps the last one, see \cite[Proposition 1.3]{2016:BassinoNicaudWeilb}.
\begin{proposition}\label{prop: basic properties ctp}
Let $\vec g$ and $\vec h$ be tuples of reduced words in $F(A)$, with the central tree property.
\begin{itemize}
\item[(1)] The Stallings graph $\Gamma(\langle \vec h\rangle)$ consists of a \emph{central tree}, which can be identified with the tree of prefixes of length $t = \lcp(\vec h)$ of the elements of $\vec h^\pm$, and of pairwise edge-disjoint paths, one for each $h_i$, from the length $t$ prefix of $h_i$ to the length $t$ prefix of $h_i\inv$.

\item[(2)] $\Gamma(\langle \vec h\rangle)$ can be computed in linear time and $\vec h$ freely generates $\langle \vec h\rangle$.

\item[(3)] If $\langle\vec h\rangle = \langle\vec g\rangle$, then $\vec h^\pm$ and $\vec g^\pm$ differ only by the order of their components.

\item[(4)] If, in addition, $3\lcp(\vec h) < \min(\vec h)$ and every word of length at most $\frac12(\min(\vec h) - 2\lcp(\vec h))$ has at most one occurrence as a factor of a word in $\vec h^\pm$, then $\langle\vec h\rangle$ is malnormal and pure.
\end{itemize}
\end{proposition}

Let us now turn to asymptotic properties. Just like we were dealing with properties of finite presentations in Section~\ref{sec: presentations} and not with properties of finitely presented groups (with the exception of 1-relator groups and $k$-relator quotients of the modular group, see Propositions~\ref{prop: statistics for 1-relator groups} and~\ref{prop: statistics for quotients of M}), we first discuss asymptotic properties of tuples of generators of a subgroup (see however Proposition~\ref{prop: statistics for k-generator subgroups} and Corollary~\ref{cor: malnormality is generic} below). We will see in Section~\ref{sec: random Stallings graphs} another way of approaching the asymptotic properties of finitely generated subgroups of $F(A)$.

In analogy with the tuples of cyclically reduced words used as relators, we can distinguish here:
\begin{itemize}
\item the density model, where at density $d$, $\P_n$ is the uniform probability law with support the $|\red_{\le n}|^d$-tuples of words in $\red_{\le n}$;
\item and the few generators model, where an integer $k\ge 1$ is fixed, and $\P_n$ is the uniform probability law with support $\red_{\le n}^k$.
\end{itemize}
Then we have the following result \cite[Propositions 3.21 and 3.22]{2016:BassinoNicaudWeilb}.

\begin{theorem}\label{thm: ctp and malnormal}
Let $0 < d < 1$.

If $d < \frac14$, then at density $d$, a tuple of reduced words $\vec h$ exponentially generically has the central tree property and in particular, it freely generates $\langle\vec h\rangle$.

If instead $d > \frac14$, then at density $d$, $\vec h$ exponentially generically does not have the central tree property.

If $d < \frac1{16}$, then at density $d$, a tuple of reduced words $\vec h$ exponentially generically generates a malnormal and pure subgroup.
\end{theorem}

It is immediate that, if every element of $\vec g$ is also an element of $\vec h$, and $\vec h$ has the central tree property, then so does $\vec g$. In that case, it is not hard to show also that $\langle\vec g\rangle$ is malnormal if $\langle\vec h\rangle$ is (see for instance \cite[Proposition 1.5]{2016:BassinoNicaudWeilb}). Then Theorem~\ref{thm: ctp and malnormal} yields the following corollary, which was already observed by Arzhantseva and Ol'shanski\u\i\ \cite{1996:ArzhantsevaOlshanskii} for the free generation statement, and Jitsukawa \cite{2002:Jitsukawa} for the malnormality statement.

\begin{corollary}\label{cor: ctp and malnormal}
In the few generators model, a tuple of reduced words exponentially generically has the central tree property, it is a basis of the subgroup it generates, and this subgroup is malnormal and pure.
\end{corollary}

We now see how to use the rigidity property in Proposition~\ref{prop: basic properties ctp} (3) to discuss asymptotic properties of subgroups themselves, and not of tuples of generators, at least in the few generators model. This is in the same spirit as in Propositions~\ref{prop: statistics for 1-relator groups} and~\ref{prop: statistics for quotients of M} above.

Fix $k \ge 1$. In the $k$-generator model for tuples, the set $S$ (in the terminology of Section~\ref{sec: models of randomness}) is $\red^k$ and $\P_n$ is the uniform probability law with support $S_n = \red_{\le n}^k$. Now consider the set $T$ of all $k$-generated subgroups of $F(A)$, the set $T_n$ of subgroups of the form $\langle\vec h\rangle$ for some $\vec h\in S_n$ and the probability law $\Q_n$ which is uniform on $T_n$. We call this the \emph{$k$-generator model for subgroups}.

\begin{proposition}\label{prop: statistics for k-generator subgroups}
Let $X$ be a property of $k$-generator subgroups of $F(A)$, that is, $X$ is a subset of $T$. Let $Y = \{\vec h\in \red^k \mid \langle\vec h\rangle \in X\}$. If $Y$ is negligible (resp. exponentially negligible) in the $k$-generator model for tuples, then so is $X$, in the $k$-generator model for subgroups. The same statement holds for genericity instead of negligibility.
\end{proposition}

\proof
Let $P$ be the set of tuples with the central tree property. By Corollary~\ref{cor: ctp and malnormal}, there exist $C,d > 0$ such that $\P_n(P) > 1-Ce^{-dn}$. Moreover, by Proposition~\ref{prop: basic properties ctp} (3), if $\vec h \in P$, there are at most $2^kk!$ elements of $P$ which generate the subgroup $\langle\vec h\rangle$.

If $Z$ is the set of subgroups of $F(A)$ of the form $\langle\vec h\rangle$ such that $\vec h\in P$, one shows as in the proof of Proposition~\ref{prop: statistics for 1-relator groups} that $\Q_n(X) \le \Q_n(X \cap Z) + \Q_n(T \setminus Z)$, and that both terms of this sum vanish exponentially fast.
\eop

The following corollary immediately follows from Corollary~\ref{cor: ctp and malnormal}.

\begin{corollary}\label{cor: malnormality is generic}
Let $k\ge 1$. In the $k$-generator model for subgroups, malnormality and purity are exponentially generic.
\end{corollary}

\begin{remark}
The proof of Proposition~\ref{prop: statistics for k-generator subgroups} does not extend to the density model: if the number of elements of a tuple $\vec h$ is a function $k(n)$ that tends to infinity, the multiplying fact $2^kk!$ is not a constant anymore, and neglibibility for $X$ is obtained only if $\P_n(Y)$ vanishes very fast (namely, if $\P_n(Y) = o(2^{k(n)}k(n)!)$).
\end{remark}

We conclude this section with a discussion of the height of the central tree of the Stallings graph of $\langle\vec h\rangle$ (that is: the parameter $\lcp(\vec h)$) for a random choice of $\vec h$. Arzhantseva and Ol'shanski\u\i\ \cite{1996:ArzhantsevaOlshanskii} showed that in the few generators model, the height of the central tree (namely the parameter $\lcp(\vec h)$) is exponentially generically at most $\alpha n$, for any $\alpha > 0$. It is in fact generically much smaller, see \cite[Proposition 3.24]{2016:BassinoNicaudWeilb}.

\begin{proposition}\label{prop: very small tree}
Let $f$ be an unbounded non-decreasing integer function and let $k\ge 1$. The following inequality holds generically for a tuple $\vec h$ chosen randomly in the $k$-generator model: $\lcp(\vec h) \le f(n)$.
\end{proposition}

This implies that, generically in the few generators model, for tuples as well as for subgroups, the proportion of vertices of $\Gamma(\langle\vec h\rangle)$ that lie in the central tree (at most $2r(2r-1)^{\lcp(\vec h) -1}$) tends to 0 (apply Proposition~\ref{prop: very small tree} with, say, $f(n) = \log\log n$).

%%%%%%%%%%%%%%%%
\subsection{Random Stallings graphs}\label{sec: random Stallings graphs}

Another point of view on random subgroups of $F(A)$ relies on the observation that each finitely generated subgroup corresponds to a unique Stallings graph, and that these graphs admit an intrinsic combinatorial characterization, as reduced rooted $A$-graphs (see Section~\ref{sec: stallings}). The problem of drawing a random subgroup can therefore be reduced to the problem of drawing a random reduced rooted $A$-graph.

When considering such graphs, it is natural to measure their size by their number of vertices (the number of edges of such a graph of size $n$ lies between $n-1$ and $2|A|n$). By extension, we say that the size of a subgroup $H$, written $|H|$, is the size of its Stallings graph $\Gamma(H)$. Then we consider the set $S$ of all Stallings graphs over alphabet $A$ (that is: of all the reduced rooted $A$-graphs), and the uniform probability law $\P_n$ with support the Stallings graphs with $n$ vertices. This is called the \emph{graph-based model} for subgroups of $F(A)$.

\paragraph{Implementation of the graph-based model}
The problem of drawing a tuple of reduced words uniformly at random is easily solved: one draws each word independently, one letter at a time, with $2r = |\tilde A|$ choices for the first letter, and $2r-1$ choices for each of the following letters.

Drawing (a tuple of) cyclically reduced words uniformly at random is also done in a simple way. Indeed, the probability that a random reduced word of length $n$ is cyclically reduced tends to $\frac{2r-1}{2r}$ when $n$ tends to infinity, and we can use a rejection algorithm: repeatedly draw a reduced word until that word is cyclically reduced. The expected number of draws tends to $\frac{2r}{2r-1} = 1 + \frac1{2r-1}$.

Drawing a Stallings graph with $n$ vertices is a less immediate task. Bassino \emph{et al.} \cite{2008:BassinoNicaudWeil} use a recursive method and the tools of analytic combinatorics to solve it in an efficient manner: they give a rejection algorithm with expected number of draws $1+o(1)$, which requires a linear time precomputation, and takes linear time for each draw. These linear time bounds are evaluated in the RAM model; in the bit complexity model, the precomputation is done in time $\O(n^2\log n)$ and each draw is done in time $\O(n^2\log^2 n)$ (see \cite[Section 3]{2008:BassinoNicaudWeil}).

\paragraph{Asymptotic properties of subgroups in the graph-based model}
The following is a combination of \cite[Section 2.4, Corollary 4.1]{2008:BassinoNicaudWeil} and \cite[Corollary 4.8 and Theorems 5.1 and 6.1]{2013:BassinoMartinoNicaud}. We say that a property $X$ is \emph{super-polynomially negligible} (resp. \emph{generic}) if $\P_n(X)$ is $\O(n^{-k})$ (resp. $1 - \O(n^{-k})$) for every positive integer $k$.
 
\begin{theorem}\label{thm: asymptotic graph based}
Let $r = |A|$.

\begin{itemize}
\item[(1)] The number of subgroups of $F(A)$ of size $n$ is asymptotically equivalent to
$$\frac{(2e)^{-r/2}}{\sqrt{2\pi}}\ e^{-(r-1)n+2r\sqrt n}\ n^{(r-1)n+\frac{r+2}4}.$$
\item[(2)] The expected rank of a size $n$ subgroup of $F(A)$ is $(r-1)n-r\sqrt n+1$, with standard deviation $o(\sqrt n)$.
 
\item[(3)] In the graph-based model, a random subgroup of $F(A)$ of size $n$ is generically neither malnormal nor pure: it is malnormal (resp. pure) with vanishing probability $\O(n^{-\frac r2})$.
 
\item[(4)] The probability that a subgroup of $F(A)$ of size $n$ avoids all the conjugates of the elements of $A$ tends to $e^{-r}$.

\item[(5)] The probability that a subgroup of $F(A)$ of size $n$ has finite index admits an $\O(n^{\frac r4} e^{-2r\sqrt n})$ upper bound. In particular, this class of subgroups is super-polynomially negligible.

\item[(6)] In the graph-based model, the quotient of $F(A)$ by the normal closure of a random subgroup is generically trivial.
\end{itemize}
 \end{theorem}
 
Some of these results call for comments, especially in comparison with the results reported in Section~\ref{sec: ctp}. As discussed at the very end of that section, in the Stallings graph of a subgroup taken at random in the few generators model, the immense majority of vertices are on the outer loops, adjacent to exactly two edges. In fact, since the rank of a subgroup is the difference between the number of edges and the number of vertices plus 1, the ratio between the number of edges and vertices tends to 1 in the few generators model (Proposition~\ref{prop: basic properties ctp}~(2) and Corollary~\ref{cor: ctp and malnormal}), and it tends to $|A|-1$ in the graph-based model (Theorem~\ref{thm: asymptotic graph based}~(2)). Observe that the minimum and maximum possible values for this ratio are $1$ and $|A|$: in intuitive terms, the Stallings graph of a random group is spase in the few generators model, and rather full in the graph based model. In other words, there are many more loops, including short loops, in the latter model, whereas in the $k$ relator model, there are only $k$ loops, and they are all very long: using close to a $\frac1k$ proportion of the edges. This is the feature that is exploited in \cite{2013:BassinoMartinoNicaud} to show that the property in Theorem~\ref{thm: asymptotic graph based}~(4) is exponentially negligible in the few generators model, and indeed in the density model at densities $d<\frac14$. Similarly, generically in the graph based model, a Stallings graph has a cycle labeled by a power of a letter, and hence the corresponding subgroup is neither malnormal nor pure (Theorem~\ref{thm: asymptotic graph based}~(3)). This is a very rough sufficient reason for a subgroup to fail being malnormal or pure, and the probability of this property may well vanish faster than stated above. A refinement of this result (namely the fact that for each letter $a$, the lengths of the cycles labeled by a power of $a$ are relatively prime) leads to Theorem~\ref{thm: asymptotic graph based}~(6). In this respect, we see that drawing uniformly at random the Stallings graph of the subgroup generated by a tuple of relators is not a fruitful avenue, to discuss 'typical' properties of finitely presentations.
 
Finally, we note that the estimates in Theorem~\ref{thm: asymptotic graph based}~(1) and~(5) can be seen as an extension of the study of subgroup growth, see in particular Lubotzky and Segal \cite{2003:LubotzkySegal}.

%%%%%%%%%%%%%%%%
\subsection{Whitehead minimality}\label{sec: Whitehead minimality}

The following property of a subgroup $H$ of $F(A)$ has already been mentioned: we say that $H$ is \emph{Whitehead minimal} (resp. \emph{strictly Whitehead minimal}) if $|\phi(H)| \ge |H|$ (resp. $|\phi(H)| > |H|$) for every non length-preserving automorphism $\phi$ of $F(A)$. This property plays an important role in the solution of the automorphic orbit problem, to decide whether two subgroups are in the same orbit under the automorphism group of $F(A)$, as shown by Gersten \cite[Corollary 2]{1984:Gersten}, in an extension of the famous Whitehead peak reduction theorem \cite{1936:Whitehead} (see also \cite[Section 1.4]{1977:LyndonSchupp}) from elements of $F(A)$ to finitely generated subgroups.

Note that a cyclic subgroup $H = \langle u\rangle$ is (strictly) Whitehead minimal if and only if the word $u$ is (strictly) Whitehead minimal in the sense discussed in Section~\ref{sec: 1-relator}. As mentioned there, Kapovich \emph{et al.} proved that strictly Whitehead minimal cyclically reduced words are exponentially generic in $F(A)$ \cite[Theorem A]{2006:KapovichSchuppShpilrain}.

This can be generalized to all finitely generated subgroups. Since the Stallings graph of a Whitehead minimal subgroup must be cyclically reduced, the graph based model must be restricted (in the natural way) to these graphs. If we consider instead the few generators model, we note that being cyclically reduced is not a generic property (see \cite[Proposition 4.6]{2016:BassinoNicaudWeil}): here too, the few generators model must be restricted to tuples of cyclically reduced words, that is, to the few relator model of Section~\ref{sec: presentations}. Under these restrictions, Bassino \emph{et al.} proved that strict Whitehead minimality is generic both in the graph based and in the few generators models \cite[Theorems 3.1 and 4.1]{2016:BassinoNicaudWeil}.

\begin{theorem}
Strict Whitehead minimality is super-polynomially generic for the uniform distribution of cyclically reduced Stallings graphs.

The same property is exponentially generic in the few relator model, restricted to tuples of cyclically reduced words.
\end{theorem}

\begin{remark}
The reasons for genericity are different for the two models, due to the very different expected geometry of a random Stallings graph: in the few generator models, it is very sparse and most of its vertices are on very long loops, whereas the graph is fuller and has many short loops in the graph-based model. See \cite{2016:BassinoNicaudWeil} for more details.
\end{remark}

%%%%%%%%%%%%%%%%
\subsection{Random subgroups of non-free groups}\label{sec: subgroups of other groups}

Let us first return to the few generators model, but for subgroups of some fixed, non-free $A$-generated group $G$. Here, the probability laws $\P_n$ we consider are the uniform probability laws with support $(\tilde A^{\le n})^k$ for some fixed $k \ge 1$: that is, we draw uniformly at random $k$-tuples of words of length at most $n$, that are not necessarily reduced.

Gilman \emph{et al.} show the following proposition \cite[Theorem 2.1]{2010:GilmanMiasnikovOsin}. Recall that a group is \emph{non-elementary hyperbolic} if it is hyperbolic and does not have a cyclic, finite index subgroup.

\begin{proposition}\label{prop: GMO}
Let $G$ be a non-elementary hyperbolic group and let $k\ge 1$. Then for any choice of generators $A$ of $G$ and any onto morphism $\pi\colon F(A)\to G$, exponentially generically in the $k$-generator model, a tuple $\vec h$ of elements of $F(A)$ is such that $\pi(\vec h)$ freely generates a free, quasi-convex subgroup of $G$.
\end{proposition}

Note that a free group $F(A)$ is non-elementary hyperbolic if $|A|\ge 2$: thus Proposition~\ref{prop: GMO} generalizes part of Corollary~\ref{cor: ctp and malnormal}, since the latter is only relative to the standard set of generators of $F(A)$.

Say that a group $G$ has the \emph{(exponentially) generic free basis property} if, for every choice of generators $A$ of $G$ and every onto morphism $\pi\colon F(A)\to G$, for every integer $k \ge 1$, the $\pi$-image of a $k$-tuple $\vec h$ of elements of $\tilde A^*$ (exponentially) generically freely generates a free subgroup of $G$ (in the $k$-generator model). Proposition~\ref{prop: GMO} states that non-elementary hyperbolic groups have the exponentially generic free basis property. Gilman \emph{et al.} \cite{2010:GilmanMiasnikovOsin} and Myasnikov and Ushakov \cite{2008:MyasnikovUshakov} note that this property is preserved as follows: if $\phi\colon G_1\to G_2$ is an onto morphism and $G_2$ has the (exponentially) generic free basis property, then so does $G_1$. For instance, non abelian right-angled Artin groups and pure braid groups $PB_n$ ($n\ge 3$) have the exponentially generic free basis property, since they admit morphisms onto a rank 2 free group (see \textit{e.g.} \cite{1974:Birman} for $PB_n$).

Proposition~\ref{prop: GMO} can be used also to show the following result \cite[Theorem 2.2]{2010:GilmanMiasnikovOsin} on the membership problem in subgroups -- a problem which is, in general, undecidable in hyperbolic groups \cite{1982:Rips}.

\begin{corollary}
Let $G$ be a non-elementary hyperbolic $A$-generated group, let $\pi$ be a surjective morphism from $A^*$ onto $G$ and let $k\ge 1$. There exists an exponentially generic set $X$ of $k$-tuples of words in $\tilde A^*$ and a cubic time algorithm which, on input a $k$-tuple $\vec h$ and an element $x\in \tilde A^*$, decides whether $\vec h \in X$, and if so, solves the membership problem for $\pi(x)$ and $\pi(\vec h)$, that is, decides whether $\pi(x) \in \langle\pi(\vec h)\rangle$.
\end{corollary}

There is no study as yet of asymptotic properties of subgroups of non free groups using a graph based model, in the spirit of Section~\ref{sec: random Stallings graphs}. Let us however mention that recent results may open the way towards such a study: Kharlampovich \emph{et al.} \cite{2016:KharlampovichMiasnikovWeil} effectively construct Stallings graphs which are uniquely associated with each quasi-convex subgroup of a geodesically automatic group, e.g. hyperbolic groups, right-angled Artin groups. Like in the free group case, this has a large number of algorithmic consequences. It may be difficult to combinatorially characterize these graphs in general, and to design random generation algorithms or to explore their asymptotic properties. But it may be possible to tackle this task for specific groups or classes of groups.

In fact, somewhat earlier results already gave more efficient and more combinatorially luminous constructions, for amalgams of finite groups (Markus-Epstein \cite{2007:Markus-Epstein}) and for virtually free groups (Silva \emph{et al.} \cite{2015:SilvaSoler-EscrivaVentura}). Note that both classes of groups are locally quasi-convex, and these constructions apply to all their finitely generated subgroups.

%%%%%%%%%%%%%%%%%%%
\section{Non-uniform distributions}\label{sec: non-uniform distributions}

In this final section, we introduce non-uniform distributions, both for relators and generators, as explored by Bassino \emph{et al.} \cite{2016:BassinoNicaudWeilb}. We keep the idea of randomly drawing tuples of words by independently drawing the elements of the tuple, but we relax the distribution on the lengths of the tuples and on the lengths of the words, and we use non-uniform probability laws of probability on each $\red_n$ (resp. $\cred_n$).

More precisely, the model of randomness is the following \cite{2016:BassinoNicaudWeilb}. For each $n\ge 0$, let $\R_n$ be a law of probability on $\red_n$ (or $\cred_n$ if we are dealing with presentations) and let $\T_n$ be a law of probability on the set of tuples of positive integers. If $\vec h = (h_1, \ldots, h_k)$ is a tuple of words, let $|\vec h| = (|h_1|, \ldots, |h_k|)$. Together, $(\R_n)_n$ and $(\T_n)_n$ define a sequence of probability laws $\P_n$ on the set of tuples of (cyclically) reduced words as follows:
$$\P_n(\vec h) = \T_n(|\vec h|) \ \prod_i \R_{|h_i|}(h_i).$$
Note that this includes the density and the few generators (relators) models discussed in Sections~\ref{sec: presentations} and~\ref{sec: random subgroups}: for instance, in the $k$-generator model, $\R_n$ is the uniform distribution on $\red_n$ and $\T_n$ is the distribution with support the $k$-tuples of integers between 0 and $n$, each with probability
$$\T_n(\ell_1,\ldots,\ell_k) = \prod_{i=1}^k\frac{|\red_{\ell_i}|}{|\red_{\le n}|}.$$

%%%%%%%%%%%%%%%%
\subsection{Prefix-heavy distributions}\label{sec: prefix-heavy}

For each word $u\in \red$, denote by $\calP(u)$ the set of reduced words starting with $u$, that is, $\calP(u) = u\tilde A^* \cap \red$. For $C \ge 1$ and $0 < \alpha < 1$, say that the sequence of probability laws $(\R_n)_n$ (each with support in $\red_n$) is \emph{prefix-heavy with parameters $(C,\alpha)$} if, for all $u,v\in \red$, we have
$$\R_n(\calP(uv) \mid \calP(u)) \enspace\le\enspace C \alpha^{|v|}.$$
This definition captures the idea that the probability of a prefix-defined set (a set of the form $\calP(u)$) decreases exponentially fast with the length of $u$. It is satisfied by the sequence of uniform probability laws on the $\red_n$ ($n\ge 0$).

If $(\P_n)_n$ is a sequence of laws of probability on tuples of reduced words, defined as above by sequences $(\R_n)_n$ and $(\T_n)_n$ of probability laws on words and on tuples of integers, and if $(\R_n)_n$ is prefix-heavy with parameters $(C,\alpha)$, then we say that $(\P_n)_n$ is \emph{prefix-heavy} as well, with the same parameters.

Under this hypothesis, Bassino \emph{et al.} obtain a series of general results \cite[Theorems 3.18, 3.19 and 3.20]{2016:BassinoNicaudWeilb}, summarized as follows. If $\vec h = (h_1, \ldots, h_k)$ is a tuple of reduced words, we let $\size(\vec h) = k$, $\min(\vec h) = \min\{|h_i| \mid 1\le i \le k\}$ and $\max(\vec h) = \max\{|h_i| \mid 1\le i \le k\}$. Let us say, also, that $(\R_n)_n$ and $(\P_n)_n$ \emph{do not ignore cyclically reduced words} if $\liminf \R_n(\cred_n) > 0$.

\begin{theorem}\label{thm: general prefix-heavy}
Let $(\P_n)_n$ be a sequence of probability laws on tuples of reduced words, which is prefix-heavy with parameters $(C,\alpha)$, with $C\ge 1$ and $0 < \alpha < 1$. Let $0 < \lambda < \frac12$.
\begin{itemize}
\item If the random variable $\size^2 \alpha^{\frac\min2}$ is increasingly small --- more precisely, if there exists a sequence $(\eta_n)_n$ tending to 0, such that $\P_n(\size^2 \alpha^{\frac\min2} > \eta_n)$ tends to 0 ---, then a random tuple of reduced words generically satisfies the central tree property, and freely generates a subgroup of $F(A)$.

\item If there exists a sequence $(\eta_n)_n$ tending to 0, such that $\P_n(\size^2 \max^2 \alpha^{\frac\min8} > \eta_n)$ tends to 0, then a random tuple of reduced words generically generates a malnormal subgroup of $F(A)$.

\item Let $0 < \lambda < \frac12$. If the sequence $(\P_n)_n$ does not ignore cyclically reduced words and if there exists a sequence $(\eta_n)_n$ tending to 0, such that $\P_n(\size^2 \max^2 \alpha^{\lambda\min} > \eta_n)$ tends to 0, then a random tuple of cyclically reduced words generically satisfies the small cancellation property $C'(\frac16)$.
\end{itemize}
In all three statements, exponential genericity is guaranteed if the vanishing sequences converge exponentially fast to 0.
\end{theorem}

The technical aspect of these statements is due to the very general nature of the random model considered. In the next section, we discuss a more specific model, where the $\R_n$ are generated by a Markovian scheme.

%%%%%%%%%%%%%%%%
\subsection{Markovian automata}\label{sec: markovian automata}

When it comes to drawing words at random, an automaton-theoretic model comes naturally to mind. Bassino \emph{et al.} introduce the following notion: a \emph{Markovian automaton} $\calA$ over a finite alphabet $X$ consists in a finite deterministic transition system $(Q, \cdot)$ (that is: an action of the free monoid $X^*$ on the finite set $Q$, or seen otherwise, a deterministic finite state automaton over alphabet $X$ without initial or terminal states), an initial probability vector $\gamma_0 \in [0,1]^Q$ (that is: $\sum_{p\in Q}\gamma_0(p) = 1$), and a stochastic matrix $M \in [0,1]^{Q\times X}$ (that is, a matrix where each column  is a probability vector) such that $M(p,x) > 0$ if and only if $p\cdot x$ is defined.

Such a scheme defines a sequence $(\R_n)_n$ of laws of probability, over each set $X^n$ ($n\ge 0$), as follows:
$$\R_n(x_1\cdots x_n) = \sum_{p\in Q} \gamma_0(p) M(p,x_1) M(p\cdot x_1,x_2) \cdots M(p\cdot (x_1\cdots x_{n-1}),x_n).$$
Note that the union over $n$ of the support sets of the $\R_n$ is always a prefix-closed rational language: that accepted by the transition system $(Q,\cdot)$, with initial states the support of $\gamma_0$ and all states final.

\begin{example}
For instance, if $Q = \tilde A$, if for each $a,b\in \tilde A$, $a\cdot b$ is defined whenever $b\ne a\inv$, and equal to $b$ when defined, if the entries of $\gamma_0$ are all equal to $\frac1{2r}$ and if the non-zero entries of $M$ are all equal to $\frac1{2r-1}$, then $\R_n$ is the uniform probability law on $\red_n$.

The Markovian automata in Figure~\ref{fig: 2 automata}
also yield the uniform probability law (at fixed length) on two languages which both provide unique %
\begin{figure}[htbp]
\begin{center}
\null\hfill
$(\calA)$\includegraphics[scale=1]{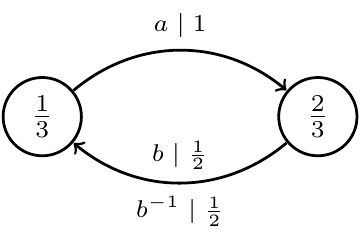}
\hfill
\includegraphics[scale=1]{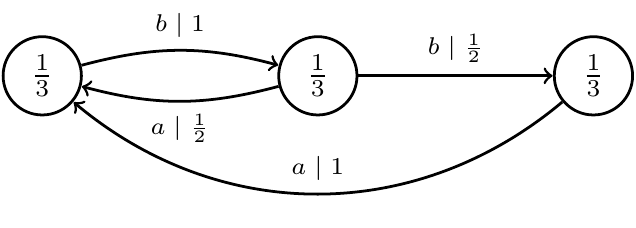}$(\calA')$
\hfill\null
\caption{\small Markovian automata $\calA$ and $\calA'$. Transitions are labeled by a letter and a probability, and each state is decorated with the corresponding initial probability.\label{fig: 2 automata}}
\end{center}
\end{figure}
representatives for the elements of the modular group (see Section~\ref{sec: rigidity}): the support of $\calA$ is the set of words over alphabet $\{a,b,b\inv\}$ without occurrences of the factors $a^2$, $b^2$, $(b\inv)^2$, $bb\inv$ and $b\inv b$ (the shortlex geodesics of the modular group), and the support of $\calA'$ consists of the words on alphabet $\{a,b\}$, without occurrences of $a^2$ or $b^3$.
\end{example}

A first set of results is obtained by specializing Theorem~\ref{thm: general prefix-heavy} to the case where the sequence $(\P_n)_n$ is induced by a Markovian automaton $\calA$. If $0 < \alpha < 1$, we introduce the $\alpha$-density model with respect to $\calA$, in analogy with Sections~\ref{sec: density model} and~\ref{sec: ctp}: at density $d < 1$, the sequence $(\P_n)_n$ is induced by the sequences $(\R_n)_n$, induced by $\calA$, and $(\T_n)_n$, where the support of $\T_n$ is reduced to the $\alpha^{dn}$-tuple $(n,\ldots,n)$. The usual density model corresponds to $\alpha = \frac1{2r-1}$.

The following is a generalization of the results in Section~\ref{sec: ctp} \cite[Proposition 4.3 and Corollary 4.5]{2016:BassinoNicaudWeilb}.

\begin{theorem}
Let $\calA$ be a Markovian automaton.

If $\calA$ does not have a cycle with probability 1, then the induced sequence of probability laws on $\red$ is prefix-heavy, with computable parameters $(C,\alpha)$.

If that is the case, then in the density model with respect to $\calA$, at $\alpha$-density $d < \frac14$, a tuple of reduced words exponentially generically has the central tree property.

And at $\alpha$-density $d < \frac1{16}$, a tuple of reduced words exponentially generically generates a malnormal subgroup.
\end{theorem}

We get more precise results if the Markovian automaton $\calA$ is ergodic, that is, if its underlying graph is strongly connected and if, for every large enough $n$, there are paths of length $n$ from every state to every other one. In that situation, it is well known that $\calA$ has a stationary vector $\tilde\gamma \in [0,1]^Q$, and we let $(\tilde\R_n)_n$ be the sequence of probability laws defined by $\calA$ with initial vector $\tilde\gamma$ instead of $\gamma_0$. We say that $\calA$ is \emph{non-degenerate} if $\sum_{a\in \tilde A}\R_n(a)\tilde\R_n(a\inv) \ne 1$. Finally we define the \emph{coincidence probability} $\alpha_{[2]}$ of $\cal A$ as follows: let $M_{[2]}$ be the $((Q \times \tilde A) \times (Q \times \tilde A))$-matrix with entries
$$M_{[2]}((p,a),(q,b)) = \begin{cases} \gamma(p,b)^2 & \textrm{if $p\cdot b = q$,} \\ 0 & \textrm{otherwise.}\end{cases}$$
Then $\alpha_{[2]}$ is the largest eigenvalue of $M_{[2]}$. Bassino \emph{et al.} proved the following phase transition result, which generalizes Theorem~\ref{thm: transition for small cancellation} and part of Theorem~\ref{thm: transition for hyperbolic}, see \cite[Propositions 4.9 and 4.14, Theorem 4.15]{2016:BassinoNicaudWeilb}.

\begin{theorem}
Let $\calA$ be a non-degenerate ergodic Markovian automaton. Then the induced sequence of probability laws on $\red$ is prefix-heavy with parameters $(C,\sqrt{\alpha_{[2]}}$), for some computable $C \ge 1$, and it does not ignore cyclically reduced words.

In particular, in the density model with respect to $\calA$, at $\alpha_{[2]}$-density $d < \frac18$, a tuple of reduced words exponentially generically has the central tree property; and at $\alpha_{[2]}$-density $d < \frac1{32}$, a tuple of reduced words exponentially generically generates a malnormal subgroup.

Moreover, let $0 < \lambda < \frac12$. Then at $\alpha_{[2]}$-density $d < \frac\lambda2$, a tuple of cyclically reduced words exponentially generically satisfies Property $C'(\lambda)$. And at $\alpha_{[2]}$-density $d > \frac\lambda2$, it exponentially generically does not satisfy Property $C'(\lambda)$.

Finally, at $\alpha_{[2]}$-density $d > \frac12$, a tuple $\vec h$ of cyclically reduced words exponentially generically presents a degenerate group, in the following sense: let $B\subseteq \tilde A$ be the set of letters which label a transition in $\calA$ and let $D = A \setminus (B \cup B\inv)$. Then $\langle A \mid \vec h\rangle$ is equal to the free group of rank $|D|+1$ if $B \cap B\inv = \emptyset$, and otherwise to $F(D)\ast\Z/\!2\Z$ if $n$ is even, $F(D)$ if $n$ is odd.
\end{theorem}

%%%%%%%%%%%%%%%%%%%%%%
%%%%%%%%%%%%%%%%%%%%%%
{\small\bibliographystyle{abbrv}

}

%%%%%%%%%%%%%%%%
%%%%%%%%%%%%%%%%
%%%%%%%%%%%%%%%%
\end{document}